# Generalization and Alternatives of Kaprekar's Routine

Florentin Smarandache
University of New Mexico, Gallup Campus, USA

**Abstract.**
We extend Kaprekar's Routine for a large class of applications. We also give particular examples of this generalization as alternatives to Kaprekar's Routine and Number. Some open questions about the length of the iterations until reaching either zero or a constant or a cycle, and about the length of the cycles are asked at the end.

1. **Generalization of Kaprekar's Routine.**

Let $f$ be an operator that maps a finite set $A = \{a_1, a_2, ..., a_p\}$, with $p \geq 1$, into itself:

$$f: A \rightarrow A.$$

Then, for any value $a \in A$ we have $f(a) \in A$ too. If we iterate this operator multiple times we get a chain like this:

$a \in A$ involves $f(a) \in A$, which involves $f(f(a)) \in A$ ($f$ iterated with itself twice), which in its turn involves $f(f(f(a))) \in A$, and so on, $f(...f(a)...) \in A$ { where, in the last case, $f$ has been iterated with itself $r \geq 1$ times; let's denote it by $f_r(a)$ }. Let's also denote, for the symmetry of notation,

$$a = f_0(a).$$

Since cardinal of $A$ is a finite positive integer, $card(A) = p < +\infty$, after at most $r = p$ iterations we get two equal iterations

$$f_i(a) = f_j(a)$$

for $i \neq j$ with $0 \leq i < j \leq r$, in the above chain $a = f_0(a), f_1(a), f_2(a), ..., f_r(a)$.

Hence, this chain of values can form a cycle:

$$f_0(a), f_1(a), ..., f_i(a), f_{i+1}(a), ..., f_{j-1}(a), f_i(a), f_{i+1}(a), ..., f_{j-1}(a),$$

with the cycle $f_i(a), f_{i+1}(a), ..., f_{j-1}(a)$ of length $j-i$.

If $j = i+1$ then the chain reaches a constant, since the cycle has only one element $f_i(a)$.

## 2. Kaprekar's Routine.

Kaprekar's Routine {see [1], [2], [3]}, extended to *k*-digit numbers, is a particular case of the above algorithm.

For k = 1 and 2 the Kaprekar's Routine reached zero, for k = 3 or 4 the algorithm reaches zero or 495, respectively zero or 6174. For k ≥ 6 it also reaches cycles.

The *set A = {0, and all k-digit positive integers}*, so its cardinal is finite, *card(A) = 90…01*, where in this number we have *k-2* zeroes.

The operator *f* does the following: arranges the digits of number *a* in descending order (*a'*) and in ascending order (*a''*) and then subtracts them: $a`-a``$. Since $a`-a``$ is also a *k*-digit number or zero (in the degenerate case when all the digits are equal), then $f(a`-a``)$ is a *k*-digit number as well or zero, therefore $f(a`-a``) \in A$. And the iteration continues in the same way. After a finite number of iterations the algorithm reaches a constant (which can be zero in the degenerate case), or a constant, or gets into a cycle.

## 3. Alternatives to Kaprekar's Routine.

3.1. Let's consider the group of permutations *P* of the digits of the *k*-digit *(k ≥ 1)* number *a*. We define the operator
$$f_P(a) = |P_1(a) - P_2(a)|$$

where $P_1$ and $P_2$ are some permutations of *k* elements *{1, 2, …, k}*, and |.| means absolute value.

And the *set A = {0; and all m-digit positive integers, m ≤ k }*.

Then the sequence of iterations reaches zero, a constant, or a cycle.

Let's see an **example**:

$P_1(\{1,2,3\}) = \{2,3,1\}$ and $P_2(\{1,2,3\}) = \{1,3,2\}$; *a* = 125, then we have |251-152| = 099;

|990-099|=891; |918-819|=099; |990-099|=891; |918-819|=099; … .

So, we reached a cycle: 125, *099, 891*, 099, 891, … .

3.2. Let's have the same group of permutations and same set *A* as in Example 3.1, but taking as operator:
$$fa_P(a) = |a - P(a)|.$$

See another example:

$P(\{1,2,3\}) = \{3,1,2\}$; $a = 125$, then $|125-512| = 387$, $|387-738| = 351$, $|351-135| = 216$,

$|216-621| = 405$, $|405-540| = 135$, $|135-513| = 378$, $|378-837| = 459$, $|459-945| = 486$,

$|486-648| = 162$, $|162-216| = 054$, $|054-405| = 351$, … .

So, we got: 125, 387, <u>351</u>, 216, 405, 135, 378, 459, 486, 162, 054, <u>351</u>, … .

3.3. Similarly if we take the operator: the absolute value of a number minus its reverse:

$$|a-reverse(a)|.$$

For example: 125, $|125-521| = 396$, $|396-693| = 297$, $|297-792| = \underline{495}$, $|495-594| = 099$,

$|099-990| = 891$, $|891-198| = 693$, $|693-396| = 297$, $|297-792| = \underline{495}$, … .

3.4. Let's consider the *Smarandache Form* of a number. A *Smarandache Palindrome* (SP) {see [4]-[13]} is a number of the form $a_1a_2...a_{n-1}a_na_na_{n-1}..a_2a_1$ or $a_1a_2... a_{n-1}a_na_{n-1}..a_2a_1$ where each $a_i$ is a positive integer of any number of digits. For example, 143431 is a SP since it can be written as 1(43)(43)1. A Smarandache Form is any number $a_1a_2...a_{n-1}a_n$, where each $a_i$ is a positive integer of any number of digits.
Then we can take the operator mapping SFs into SFs in the following way.
Example: Consider the following SF: 3-digit numbers under the SF of 1-digit and 2-digit groups: for example 5(76). Then we switch the groups and add them. Take modulo 1000 of the result.
Start with 5(76), 5(76) + (76)5 = 1342 whole module 1000 is 342 = 3(42);
3(42) + (42)5 = 767 = 7(67); 7(67) + (67)7 = 1444 whose modulo 1000 is 444 = 4(44);
4(44) + (44)4 = 888 = 8(88); 8(88) + (88)8 = 1776 whose module 1000 is 776 = 7(76);
7(76) + (76)7 = 1543 whose modulo 1000 is 543 = 5(43);
5(43) + (43)5 = 9(78); 9(78) + (78)9 = 1767 whose modulo 1000 is 767 = 7(67);
7(67) + (67)7 = 1444 whose modulo 1000 is 4(44).
We got: 5(76), 3(42), <u>4(44)</u>, 7(76), 5(43), 7(67), <u>4(44)</u>, … .

3.5. Or one consider another operator that subtracts two *k*-digit numbers in the following way: adding 1 to each digit less than 9, then subtracting 1 from each non-zero digit, then subtracting the numbers.
Example: 495, 596-384 = 212, 323 − 101 = <u>222</u>, 333-111 = <u>222</u>, … . We reached a constant.
Or add 2 to each digit strictly less than 8, and subtract 3 from each digit strictly greater than 2. Etc.

3.6. Consider any function *f* defined on the set of *k*-digit numbers whose range is the set of positive integers, and then calculate *modulo $10^k$* of the result.

There are infinitely many such operators in order to choose from nice examples.

**4. Open Questions:**
1. What is the longest number of iterations until one reaches either zero or a constant or a cycle that one can have for each case of the above generalizations?
2. What is the longest cycle that one can have for each particular case of the above generalizations?
3. In what conditions one reaches a constant, not a cycle? By cycle we understand a sequence of two or more numbers that repeat indefinitely.
4. Study the cases when *f(a) = a* for interesting particular cases of this generalization.
5. Study the case when *f(a) = 0* for interesting particular cases of this generalization.
6. If the operator defined in the above Generalization of the Kaprekar's Routine is a <u>random operator</u> (i.e. for a given *k*-digit number *a* one randomly generate another k-digit number *b*), is it still possible to reach a constant or a cycle?
   It is possible for sure to generate two equal *k*-digit numbers, *c*, after at most $10^k$ random operations, but then the next *k*-digit number following the first *c* would not necessarily be the same as the previous *k*-digit number following the previous *c*.